\input amstex
\documentstyle{amsppt}
\magnification 1000
\pagewidth{5.5 true in}
\pageheight{7.5 true in}
\hcorrection{.5 true in}
\vcorrection{.5 true in}
\overfullrule=0pt
\parskip=3pt
\addto\tenpoint{\normalbaselineskip17pt\normalbaselines}
\nologo
\def\e{Ext}
\def\h{Hom}
\def\t{Tor}
\def\c{\checkmark{}} 

\def\l{\lambda}

\topmatter
\title{Derived Functors and Hilbert Polynomials}\endtitle
\author{Emanoil Theodorescu}\endauthor
\address{Department of Mathematics, University of Kansas, Lawrence, KS, 66045}
\endaddress
\email{theodore\@math.ukans.edu}\endemail
\abstract{Let $R$ be a commutative Noetherian ring, $I$ an ideal, 
$M$ and $N$ finitely generated $R$-modules. Assume 
$V(I)\cap Supp(M)\cap Supp(N)$ consists of finitely many maximal ideals 
and let $\l (\e^i(N/I^nN,M))$ denote the length of $\e^i(N/I^nN,M)$. 
It is shown that $\l (\e^i(N/I^nN,M))$  agrees with a polynomial in $n$ 
for $n>>0$, and an upper bound for its degree is given. On the other hand, 
a simple example shows that some special assumption such as the support 
condition above is necessary in order to conclude that polynomial 
growth holds.}
\endabstract
\endtopmatter

\document
\head{Introduction}
\endhead
\quad Throughout this paper, $R$ is a commutative 
Noetherian ring, $I$ an ideal, $M$, $N$ are finitely generated $R$-modules, 
unless otherwise stated.
 It is proven in classical multiplicity theory (see \cite{BH}) that 
$\l(M/I^nM)$ is given by a polynomial in $n$ for $n\gg 0$, where $M$ is a 
finitely generated module over a local Noetherian ring $R$, $I$ an ideal 
and $\l(M/IM)$ is finite. In his paper \cite{Kod}, Vijay Kodiyalam proves 
a generalized version of this result, by showing that $\l(\t^i(M,N/I^nN))$ 
is eventually polynomial, provided $\l(M\otimes N)$ is finite, where $M$, $N$ 
are finitely generated modules over $R$. His approach, which makes 
consistent use of graded module theory, also yields a similar result in the 
case of $\e^i(M,N/I^nN)$, but fails to work in the case of 
$\e^i(N/I^nN,M)$.
  This latter case is partially treated by D. Kirby in \cite{Kir}, for $N=R$ 
and $i=i_0$, the first nonvanishing $\e$. Kirby shows that polynomial growth
follows for $\l(\e^{i_0}(R/I^n,M))$ from the mere assumption that 
$\e^{i_0}(R/I,M)$ has finite length. 
 In the present paper we show that if $M$, $N$ are finitely generated modules 
such that $V(I)\cap Supp(M)\cap Supp(N)$ is a finite set of maximal ideals 
of $R$, then, for all $i\geq 0$, $\l (Ext^i(N/I^nN,M)$ has polynomial growth 
for $n \gg 0$ (see Corollary 6). We also give upper bounds for the degrees
of the polynomials associated to the Hilbert functions occurring in
this paper.
The support condition above is a reasonable one, and generalizes 
Kodiyalam's condition $\l (N\otimes M)< \infty$ to cover cases where 
$N\otimes M$ does not have finite length and, in some situations, for $M$ 
infinitely generated. As it turns out, polynomial growth no longer holds true
for $\e^i(N/I^nN,M)$ without this kind of hypothesis, in other words, it 
does not follow from the simpler assumption that $\e^i(N/I^nN,M)$ has finite
length for all large $n$.
The paper ends with an example that shows that the condition on supports
can't be simply removed if we want to conclude that polynomial growth
holds.

\subhead{Acknowledgement}
\endsubhead
I would particularly like to thank my advisor, Prof. D. Katz, for his
guidance in writing this paper.
 
\head{Preliminaries}
\endhead

 In this section we develop the method used to establish that polynomial 
growth holds for each of the following: $\e^i(M,N/I^nN)$, $\t^i(M,N/I^nN)$ 
and $\e^i(N/I^nN,M)$.
A crucial (though not difficult) result for the method of this paper is 
Proposition 3, which is used to obtain an explicit formula for 
$\e^i(M,N/I^nN)$, $\t^i(M,N/I^nN)$ and the Matlis dual of 
$\e^i(N/I^nN,M)$, and to give a degree estimate for the polynomials
associated to these modules.
As a corollary, most of Kodiyalam's results are recovered with weaker 
hypotheses.
 Since Matlis duality is being consistently used throughout this paper, a 
Matlis duality type of lemma for complexes is proven first (for the sake of
completeness, since this type of result is well-known, at least for modules).

\proclaim{Lemma 1} Let $I$ be an ideal of $R$ and 
$\Cal C$ a complex of $R$-modules (not necessarily finitely generated),

$$
\Cal C: \quad \cdots @>>> C_{i+1} @>{\partial}_{i+1}>> C_{i} @>{\partial}_{i}>>
C_{i-1}@>>> \cdots .
$$

\noindent
Let $E$ be an injective $R$-module and write $M^{\c}$ for $Hom(M,E)$,
for any $R$-module $M$.
Then we have an isomorphism of complexs

$$
Hom(R/I,\Cal C)^{\c} \cong \ \Cal C^{\c} \otimes R/I.
$$

\endproclaim

\demo{Proof} First note that if $\Cal C$ is a complex, then 
$\Cal C \otimes R/I$ is isomorphic to the complex 

$$
\Cal C/I\Cal C: \quad \cdots @>>> C_{i+1}/IC_{i+1} 
@>{\overline{\partial}}_{i+1}>> C_{i}/IC_i @>{\overline{\partial}}_{i}>>
C_{i-1}/IC_{i-1} @>>> \cdots ,
$$

\noindent
where $\overline{\partial _i}$ are the induced maps. Next, we prove that 
if we identify $Hom(R/I,M)$ with $(0 :_M I)$, then the dual with respect 
to $E$ of the inclusion $(0 :_M I) \hookrightarrow M$ is a surjective 
map $p$ with kernel $IM^{\c}$.  

Assume $I = (a_1, \ldots,a_d)$, $a_i \in R$ and apply 
\ ' \ $\c $ \ ' \  to the exact sequence

$$
0 @>>> (0 :_M I) @>>> M @>{\pmatrix a_1\\ \vdots \\ a_d\\ \endpmatrix}>> 
\oplus _{i=1}^d M.
$$

\noindent
We get 

$$
\oplus _{i=1}^d M^{\c } @>{\pmatrix a_1 & \cdots & a_d\\ \endpmatrix}>>
M^{\c } @>p>> (0 :_M I)^{\c } @>>> 0,
$$

\noindent
which shows that the kernel of $p$ is $IM^{\c}$. This means 
that there is an isomorphism \newline
$\varphi : M^{\c}/IM^{\c} @>>> Hom(R/I,M)^{\c}$ 
fitting in the following diagram  

$$
\CD
M^{\c} @>{\pi}>> M^{\c}/IM^{\c}\\
@VpVV                           \\
Hom(R/I,M)^{\c}                {}\\
\endCD
$$

\noindent
where 
$$\varphi := p \circ \pi^{-1}\tag1$$ 
\newline is well-defined.
Now, if  $M_1@>\partial>>M_2$ is a map of $R$-modules, we get the diagram

$$
\CD
0 @>>> Hom(R/I,M_1) @>{i_1}>> M_1 \\
@.     @VV{\delta}V   @VV{\partial}V \\ 
0 @>>> Hom(R/I,M_2) @>{i_2}>>  M_2 \\
\endCD
\tag2
$$

\noindent
from which, dualizing, we get

$$
\CD
M_1^{\c}@>{p_1}>>Hom(R/I,M_1)^{\c}@>>>0\\
@AA{\partial^{\c}}A  @AA{\delta^{\c}}A @.\\
M_2^{\c}@>{p_2}>>Hom(R/I,M_2)^{\c}@>>>0.\\
\endCD
\tag3
$$

\noindent
We also have the commutative diagram

$$
\CD
M_1^{\c} @>{\pi_1}>>M_1^{\c}/IM_1^{\c} @>>>0\\
@AA{\partial ^{\c}}A    @AA{D}A @.\\
M_2^{\c} @>{\pi_2}>> M_2^{\c}/IM_2 ^{\c} @>>>0.\\
\endCD
\tag4
$$

\noindent
The fact that $\delta$ in (2) is the restriction of $\partial$ to the
corresponding submodules is equivalent to the commutativity of (2):
$\quad \partial \circ i_1=i_2 \circ \delta$.
Dualizing, we get the commutativity of (3):
$\quad  p_1 \circ \partial^{\c}=\delta^{\c} \circ p_2$.
But this says that

$$
\quad \delta^{\c}:=p_1\circ \partial^{\c}\circ p_2^{-1}\tag5
$$

\noindent
is well-defined.
Similarly, the map in (4)

$$
\quad D :=\pi_1 \circ \partial^{\c} \circ \pi_2^{-1}\tag6
$$

\noindent
is well-defined.
An easy computation now shows that
$\varphi_1 \circ D{\buildrel (1),(6) \over =}(p_1\circ \pi_1^{-1})
\circ(\pi_1\circ \partial^{\c}\circ \pi_2^{-1})
=p_1\circ \partial^{\c} \circ \pi_2^{-1} {\buildrel 
(5) \over =} \delta ^{\c} \circ p_2 \circ \pi _2 ^{-1} {\buildrel (1) 
\over =} \delta ^{\c} \circ \varphi _2$.
Applying all this to the complex $\Cal C$, with $M_1=C_i$, $M_2=C_{i-1}$,
and $\partial = \partial_i$, gives that the two complexes 
$Hom(R/I, \Cal C)^{\c}$ and $\Cal C ^{\c}/I\Cal C ^{\c}$ are isomorphic,
the isomorphism being given by $\{\varphi_i \}_i$.
\enddemo

 The following Lemma gives the required polynomial growth conclusion
in all the situations occurring in this paper.

\proclaim{Lemma 2} Let $R$ be a Noetherian ring, $I$ an ideal, $T$ an 
$R$-module, $U$, $V$, $W$ and $Z$ finitely generated submodules of $T$. 

\roster
\item"(a)"
Assume that one of the following containments hold:
$Z \subset U$, or $W \subset V$.
Then there exist finite submodules $A$, $B \subset T$ and $k \in \Bbb{N}$,
such that, for all $n \geq k$, the following intersection formula holds:
$$
\quad (U + I^nV)\cap (Z+I^nW)= A + I^{n-k}B.
$$

\item"(b)"
Assume $Z \subset U$ and $W \subset V$.
Set

$$
L_n:=\frac{U+I^nV}{Z+I^nW},\ n \gg 0,
$$ 

\noindent
and assume that $L_n$ has finite length  for all $n \gg 0$.
Then the length of $L_n$ has polynomial growth with respect to $n$, $n\gg 0$.
\item"(c)"
Assume that $R$ is local with maximal ideal $m$, and 
$J\subseteq ann_{R} (V+Z)/Z$. 
Then, with the notations and assumptions in $(b)$, 
$$
\deg \l (L_n) \leq \max \{\dim U/Z, \ell_{R/J}(I)-1\},
$$ 
\noindent
where $\ell_{R/J}(I)$ denotes the analytic spread of $I$ on the module 
$R/J$. 
\endroster
\endproclaim

\demo{Proof}
(a) Assume first that $Z \subset U$. Using the modularity lemma, we obtain
$$
(U+I^nV)\cap (Z+I^nW)=(U+I^nV)\cap I^nW + Z.
$$
Note that over the Rees ring ${\Cal R}_I:= \oplus_{n=0}^{\infty}I^n$,
$\oplus_{n=0}^{\infty}(U+I^nV)$ and $\oplus_{n=0}^{\infty}I^nW$ are graded
${\Cal R}_I$ submodules of $\oplus_{n=0}^{\infty}T$. Then so is their 
intersection, $\oplus_{n=0}^{\infty}(U+I^nV)\cap I^nW$, which is 
finite over ${\Cal R}_I$, since $\oplus_{n=0}^{\infty}I^nW$ is finite over
the Rees ring.
But then
$$
(U+I^nV)\cap I^nW =I^{n-k}((U+I^kV)\cap I^kW)
$$
for some $k \in \Bbb N$ and $n \geq k$, by the graded Artin-Rees lemma. 
The conclusion follows with $A=Z$ and $B=(U+I^kV) \cap I^kW$. 

Now assume $W \subset V$. By the modularity lemma again, we get
$(U+I^nV)\cap (Z+I^nW)=(U+I^nV)\cap Z + I^nW$, which reduces the problem to 
the case $W=0$. In order to compute $(U+I^nV)\cap Z$, consider first that
$U\cap Z=0$. Let $x \in (U+I^nV)\cap Z$. Then $x=u+\sum_r i_rv_r=z$, where
$u \in U$, $z \in Z$, $i_r \in I^n$, $v_r \in V$. We obtain 
$y:=-u+z=\sum_r i_rv_r$, $y \in (U\oplus Z)\cap I^nV=
I^{n-k}((U\oplus Z)\cap I^kV)$, for some $k$ and $n\geq k$.
Let $\pi_2 : U\oplus Z @>>> Z$ be the second projection.
Thus, $x=z=\pi_2(y)$, $x \in I^{n-k}\pi_2((U\oplus Z)\cap I^kV)$.
Conversely, let $x \in I^{n-k}\pi_2((U\oplus Z)\cap I^kV)=
\pi_2((U\oplus Z)\cap I^nV)$. This means 
$x=z$, $z \in Z$, and that for some $u \in U$, $-u+z=\sum_r i_rv_r$, 
$i_r \in I^n$, $v_r \in V$. Hence, $x=z=u+ \sum_r i_rv_r$, giving 
$x \in (U+I^nV)\cap Z$. We just proved that, if $U\cap Z=0$, then
$(U+I^nV)\cap Z=I^{n-k}B$, where $B=\pi_2((U\oplus Z)\cap I^kV)$, $\pi_2$ as
above. Finally, in the general case, denote by $\overline U$, $\overline V$,
$\overline Z$, the images in $T/(U\cap Z)$ of $U$, $V$, $Z$, respectively.
In $T/(U \cap Z)$ we have $\overline U \cap \overline Z=0$, so
$(\overline U + I^n\overline V)\cap \overline Z=I^{n-k}\overline L$,
for some $\overline L \subset T/(U \cap Z)$. If we take $L$ to be the inverse
image of $\overline L$ in $T$, we finally obtain 
$(U+I^nV)\cap Z= U\cap Z+ I^{n-k}L$, which is what we wanted.
\newline
(b) First reduce to the case $Z=0$.
Since $Z \subset U$, we have

$$
L_n\cong \frac{(U+I^nV)/Z}{(Z+I^nW)/Z} \cong \frac{U/Z+I^n(V+Z)/Z}{I^n(W+Z)/Z}.
$$

\noindent
By replacing $T$ by $T/Z$, $U$ by $(U+Z)/Z$, $V$ by $(V+Z)/Z$, 
$W$ by $(W+Z)/Z$, the inclusion $W\subset V$ still holds. 
Note now that both

$$
\frac{U+I^nW}{I^nW} \quad \text{and} \quad \frac{I^nV}{I^nW}
$$ 

\noindent
have finite length, as submodules of $L_n$.

We have an obvious short exact sequence

$$
0 @>>> \frac{U+I^nW}{I^nW} \cap \frac{I^nV}{I^nW} @>>> \frac{U+I^nW}
{I^nW} \oplus \frac{I^nV}{I^nW} @>>> L_n @>>>0,
\tag 7
$$

\noindent
which can be rewritten as

$$
0@>>> \frac{(U+I^nW) \cap I^nV}{I^nW} @>>> \frac{U+I^nW}{I^nW} 
\oplus \frac{I^nV}{I^nW} @>>> L_n @>>>0.
$$

\noindent

Thus, it suffices to check that the lengths of the first two modules 
in the above short exact sequence are given by polynomials for $n \gg 0$.  
The desired conclusion follows for $\l(I^nV/I^nW)$ since
$\oplus_{n=0}^{\infty} I^nV/I^nW$ is a finitely generated 
graded module over the Rees ring ${\Cal R}_I$.
Similarly, by using the modularity lemma and the Artin-Rees lemma, we get

$$
\frac{(U+I^nW) \cap I^nV}{I^nW}=\frac{(U\cap I^nV)+I^nW}{I^nW} \cong
\frac{U\cap I^nV}{U\cap I^nW}=\frac{I^{n-k}(U\cap I^kV)}
{I^{n-k}(U\cap I^kW)}
\tag 8
$$
 
\noindent
for some $k$ and all $n\geq k$, and so again 
$$
\oplus_{n=0}^{\infty} \frac{(U+I^nW) \cap I^nV}{I^nW}
$$

\noindent
is a graded, finitely generated module over ${\Cal R}_I$,
so the needed conclusion follows again.

Finally, using Artin-Rees, we get that for some $k$ and all $n\geq k$,

$$
\frac{U+I^nW}{I^nW} \cong \frac{U}{U\cap I^nW}= \frac{U}{I^{n-k}(U 
\cap I^kW)},
$$

\noindent
and thus

$$
\l(\frac{U+I^nW}{I^nW})=\l(\frac{U}{U\cap I^kW})+
\l(\frac{U\cap I^kW}{I^{n-k}(U\cap I^kW)}).
\tag 9
$$

\noindent
Since the last term above is known to coincide with a  polynomial
for all large $n$, by general theory, it follows that $\l(L_n)$ has polynomial
growth, as stated. 
\newline
(c) Denote the images of $U$, $V$, $W$ mod $Z$ again by $U$, $V$, $W$, 
respectively. Recalling that in the proof of (b) we have reduced modulo $Z$,
the short exact sequence in (7) shows that 
$$
\deg \l(L_n)\leq \max \{\deg \l((U+I^nW)/I^nW),\deg \l(I^nV/I^nW)\}.
$$
If $U\cap I^kW=0$ in (9), then $U$ has finite length and 
$\deg \l((U+I^nW)/I^nW)= \dim U$. If we assume 
now that $U\cap I^kV \ne 0$, then $\dim (U\cap I^kW) \leq \dim U$. 
Actually, equality holds here, since $U/(U\cap I^kW)$ has finite length, 
and thus every prime ideal in the support of $U$ must also belong to the 
support of $U\cap I^kW$. So, again, (9) gives $\deg \l((U+I^nW)/I^nW)=
\dim U$.
On the other hand, $\Cal M:=\oplus_{n=0}^{\infty} I^nV/I^nW$ is a finite 
graded module over $\Cal{\overline R}_{\overline I}$, where 
\ ` \ $\bar{\ }$   \ ' \ denotes reduction mod $J$.  
Note that, after possibly dropping a few terms of $\Cal M$, we may 
assume that $J_0:=ann_{R} I^nV/I^nW$ is an $m$-primary ideal independent of 
$n$. Then $\Cal M$ is a graded module over $\Cal {\overline R}_{\bar I}/ 
\bar {J_0}\Cal {\overline R}_{\bar I}$. Therefore $\deg \l(I^nV/I^nW)=
\dim \Cal{M}-1\leq \dim (\Cal {\overline R}_{\bar I}/ 
\bar {J_0}\Cal {\overline R}_{\bar I})-1$. Note now that 
$\dim (\Cal {\overline R}_{\bar I}/ \bar {J_0}\Cal {\overline R}_{\bar I})
=\dim (\Cal {\overline R}_{\bar I}/ \overline {m}\Cal {\overline R}_{\bar I})$,
since $J_0$ is $m$-primary, and finally this equals 
$\dim (\oplus I^n{\overline R}/mI^n{\overline R})$, the analytic spread of 
$I$ on $R/J$.
The proof of the Lemma is now complete.
\enddemo
Proposition 3 below shows that polynomial growth holds eventually, for 
a particular homology module. It actually shows that this homology  module
has the form appearing in Lemma 2(b). This will turn out to be the kind of 
formula that holds for each of the following: $\t(M,N/I^nN)$, $\e(M,N/I^nN)$,
and the Matlis dual of $\e(N/I^nN,M)$, the latter subject to  a mild 
restriction on supports. It also gives a slightly refined estimate of the 
degree of the Hilbert polynomials associated to these modules.

\proclaim{Proposition 3} Let $R$ be Noetherian and $I$ an ideal of $R$, 
$N$ a finite $R$-module, and moreover, let 

$$
\Cal{C}: \quad F_2 @>\psi>> F_1 @>\phi>> F_0
$$

\noindent
be a complex of $R$-modules.

\roster
\item"(a)"  
Assume that $F_0, F_1$ are finitely generated, and that
$H_1(\Cal C \otimes \frac{N}{I^nN})$ has finite length for 
$n \gg 0$.
Then $\l(H_1(\Cal C \otimes \frac{N}{I^nN}))$ is given by a 
polynomial for $n \gg 0$.
\item"(b)"
Assume that $F_0$, $F_1$, $F_2$ are finitely generated flat modules, and that
$H_1(\Cal C \otimes \frac{U+I^nV}{I^nW})$ has finite length for 
$n \gg 0$, where $U$, $W$ are submodules of a finitely generated 
$R$-module $V$.
Then $\l(H_1(\Cal C \otimes \frac{U+I^nV}{I^nW}))$ is given by a 
polynomial for $n \gg 0$.
\item"(c)" 
Assume that $R$ is local and the hypotheses in (a) hold. 
Then 
$$
\deg \l(H_1(\Cal{C}\otimes \frac{N}{I^nN}))\leq \max 
\{ \dim H_1(\Cal{C}\otimes N), \ell_N(I)-1\}.
$$
\noindent
If $\dim H_1(\Cal{C}\otimes N)\geq 
\ell_N(I)$, then the inequality above becomes an equality.
\endroster
\endproclaim

\demo{Proof}
(a) We may reduce to the case $N=R$, by replacing $\Cal C$ by 
$\Cal C':=\Cal C\otimes N$.\newline
Let $K:= ker \phi$, \ $ L:= im \psi$, \ so $L\subset K \subset F_1$ 
are finitely generated modules.
For all $n \geq 1$, we have the induced complexes
$$
\Cal {C'}_n: \quad \frac{F_2}{I^nF_2} @>\psi_n>> \frac{F_1}{I^nF_1} @>\phi_n>>
\frac{F_0}{I^nF_0}.
$$
\noindent
Clearly,
$$
\quad im \psi_n =\frac{L+I^nF_1}{I^nF_1}.
$$
\noindent
On the other hand,
$$
\quad ker \phi_n =\{\overline {x} \in \frac{F_1}{I^nF_1}\ : \ \phi (x) 
\in I^nF_0\}.
$$
But
$$
\phi (x) \in I^nF_0 \Longleftrightarrow \phi (x) \in I^nF_0 \cap im \phi =
I^{n-k}(I^kF_0\cap im \phi ),
$$
\noindent
by the Artin-Rees Lemma, for some $k$ and all $n \geq k$.
This means 
$$
\overline{x} \in ker \phi_n \Longleftrightarrow \phi (x)=\sum_j{i_j\phi (y_j)},
$$
\noindent
with $i_j \in I^{n-k}, y_j \in F_1, \phi (y_j) \in I^kF_0$,
so

$$
\quad \phi (x- \sum_j {i_jy_j})=0,
$$
\noindent
where
$$
y_j \in \phi^{-1} (I^kF_0)=:\tilde K,
$$
\noindent
a finitely generated submodule of $F_1$ containing $K$.
Thus, 
$\overline{x} \in ker \phi_n \Longleftrightarrow x \in K+I^{n-k}\tilde K$,
hence
$$
\quad ker \phi_n = \frac{K+I^{n-k}\tilde K}{I^nF_1}.
$$
\noindent
It follows that 

$$
H_1(\Cal C'\otimes R/I^n)= \frac{ker \phi_n}{im \psi_n} \cong \frac{K+I^{n-k}
\tilde K}{L+I^nF_1}=\frac{K+I^{n-k}\tilde K}{L+I^{n-k}I^kF_1}.
$$
\noindent
Note that $I^kF_1 \subset \tilde K$, since $\phi (I^kF_1) \subset I^kF_0$.
By taking

$$
U:=K,\ V:=\tilde K,\ Z:=L,\ W:=I^kF_1
$$
\noindent
the desired conclusion now follows from Lemma 2(b). Also note that
$U/Z=K/L=H_1(\Cal C')=H_1(\Cal C \otimes N)$.
\newline
(b)
\quad We will reduce the proof to the case where $W=V$.
Assume that we already have the result when $W=V$.
Consider the short exact sequence

$$
   0@>>>A_n@>>>B_n@>>>C_n@>>>0,
$$
where $A_n=I^nV/I^nW$, $B_n=(U+I^nV)/I^nW$, $C_n=(U+I^nV)/I^nV$, and the
maps are the canonical ones. By tensoring this with $\Cal {C}$, using the 
flatness of its modules, we get an exact sequence of complexes

$$
  0@>>>A_n\otimes \Cal{C}@>>>B_n\otimes \Cal {C}@>>>C_n\otimes \Cal {C}@>>>0.
$$
From the long exact sequence in homology, we obtain

$$
  ...@>>>H_1(A_n\otimes \Cal {C})@>\alpha_n>>H_1(B_n\otimes \Cal {C})
@>\beta_n>>H_1(C_n\otimes \Cal {C})@>\gamma_n>>H_0(A_n\otimes \Cal {C})@>>>...
$$

Moreover, $\oplus_{n=0}^{\infty} H_1(A_n\otimes \Cal {C})$ is a finite graded 
module over the Rees ring ${\Cal R}_I$ (because $\oplus_{n=0}^{\infty}A_n$ 
is so), and thus, the images of $\alpha_n$ and $\gamma_n$ are $n$-th degree 
pieces of finite ${\Cal R}_I$-modules. It follows that their lengths are  
eventually given by polynomials. The conclusion follows by exactness, since we 
assumed we have the result in the case $W=V$. In order to complete the proof, 
it remains to treat the case $W=V$.

 Denote $V\otimes F_i$, $U\otimes F_i$, by 
$\tilde{V_i}$, $\tilde{U_i}$, respectively, where 
$i \in \{0, 1, 2\}$. 
By the flatness of the $F_i$ and the right exactness of the tensor product, 
we can view $\tilde{U_i}$, $I^n\tilde{V_i}$, as submodules of $\tilde{V_i}$. 
Tensoring the complex $\Cal C$ with $(U+I^nV)/I^nV$ gives
$$
\tilde{\Cal C}_n: \quad \frac{\tilde{U_2} + I^n\tilde{V_2}}{I^n\tilde{V_2}}
@>\overline{\psi}>>\quad \frac{\tilde{U_1} + I^n\tilde{V_1}}{I^n\tilde{V_1}}
@>\overline{\phi}>>\quad \frac{\tilde{U_0} + I^n\tilde{V_0}}{I^n\tilde{V_0}},
$$
where $\overline{\phi}$ and $\overline{\psi}$ denote maps induced by 
$\tilde{\phi}:=\phi \otimes 1_V$ and $\tilde{\psi}:=\psi \otimes 1_V$, 
respectively.
Note that $\tilde{\phi}$ maps $\tilde{U_1}$ into $\tilde{U_0}$. Similarly, 
$\tilde{\psi}$ maps $\tilde{U_2}$ into $\tilde{U_1}$.
 An easy computation shows that
$$
Im \overline{\psi}=\frac{\tilde{\psi}(\tilde{U_2})+I^n\tilde{\psi}
(\tilde{V_2})+I^n\tilde{V_1}}{I^n\tilde{V_1}}=
\frac{\tilde{\psi}(\tilde{U_2})+I^n\tilde{V_1}}{I^n\tilde{V_1}} .
$$ 
We also get that
$$
Ker \overline{\phi}=\frac{\tilde{U_1}+I^n\tilde{V_1}}{I^n\tilde{V_1}}\cap 
\frac{Ker \tilde{\phi}+I^{n-l}\tilde{\phi}^{-1}(I^l\tilde{V_0})}
{I^n\tilde{V_1}},
$$
for some $l$ and all $n\geq l$.
Since it is clear that $I^l\tilde{V_1} \subset 
\tilde{\phi}^{-1}(I^l\tilde{V_0})$, by using Lemma 2(a), we obtain that 
there is a finite submodule $L$ of $\tilde{V_1}$, and some $k$, such 
that, for all $n \geq k$, we have 
$$
\quad H_1(\tilde{\Cal C})=
\frac{(\tilde{U_1}+I^n\tilde{V_1})\cap (Ker\tilde{\phi} + I^{n-k}
\tilde{\phi}^{-1}(I^k\tilde{V_0}))}
{\tilde{\psi}(\tilde{U_2})+I^n \tilde{V_1}}=
\frac{\tilde{U_1}\cap Ker\tilde{\phi}+I^{n-k}(I^k\tilde{V_1}+L)}
{\tilde{\psi}(\tilde{U_2})+I^n \tilde{V_1}}
$$
The hypotheses of Lemma 2(b) are easily seen to be satisfied with $n-k$ 
replacing $n$, since it is immediate that $\tilde{\psi}(\tilde{U_2})
\subset \tilde{U_1}\cap Ker\tilde{\phi}$.

Thus, $\l(H_1(\tilde{\Cal C}_n))$ is given by a polynomial for $n\gg 0$.
\newline
(c) Let $J:= ann_{R} N$. The proof of (a) shows that 
$H_1(\Cal{C}\otimes N/I^nN)$ has the form given in 2(b) for $L_n$.
Also note that $J$ satisfies the hypothesis in 2(c). $J$ kills every module
in $\Cal {C} \otimes N$, so it kills every subquotient of such a module.
In particular, it kills $(V+Z)/Z$, where $V$, $Z$ were defined at the end 
of 3(a). 
Then, by Lemma 2 (b), (c), and the fact that now $U/Z=H_1(\Cal C\otimes N)$,
we already know that $\deg \l(H_1(\Cal{C}\otimes N/I^nN))\leq \max 
\{\dim (H_1(\Cal{C}\otimes N)), \ell_{R/J}(I)-1\}$. It is enough to check that
$\ell_{R/J}(I)=\ell_{N}(I)$. 
By definition, $\ell_{N}(I)$ is the dimension of the graded module 
$\Cal{N}:=\oplus_{n=0}^{\infty} I^nN/mI^nN$ over the Rees ring ${\Cal R}_I$, 
or actually over $\oplus_{n=0}^{\infty} I^n(R/J)/mI^n(R/J)$, which is a 
homomorphic image of ${\Cal R}_I$. Since 
$\oplus_{n=0}^{\infty} I^n(R/J)/mI^n(R/J)\cong 
\oplus_{n=0}^{\infty} I^n/(I^n\cap J + mI^n)$, we see that 
$\oplus_{n=0}^{\infty} (I^n\cap J + mI^n)\subseteqq ann_{\Cal{R}_I}\Cal{N}$.
We are done if we prove that $ann_{\Cal{R}_I}\Cal{N}$ is, up to radical, 
equal to $\oplus_{n=0}^{\infty} (I^n\cap J + mI^n)$. Note that 
$ann_{\Cal{R}_I}\Cal{N}$ is a homogeneous ideal of $\Cal{R}_I$, and pick
$xt^k\in (It)^k$, such that $xI^nN\in mI^{n+k}N$, for all $n$
($t$ is just a ``place keeper'' in $\Cal{R}_I$).
By the determinant trick, we get $x\in \overline{mI^k}$ modulo $J$
(here the bar denotes the integral closure in $R$ of $mI^k$). Let 
$$
   x^r+a_1x^{r-1}+\cdots +a_r=j, \, a_i \in (mI^k)^i, j \in J
$$
be some integral dependence relation for $x$. Multiplying through by $t^{kr}$,
we obtain that \newline 
$(xt^k)^r\in (mI^{kr} + J\cap I^{kr})t^{kr}$, which proves what we wanted.

Assume now that $\dim(H_1(\Cal{C}\otimes N))\geq\ell_{N}(I)$. Note that
$H_1(\Cal C\otimes N)=U$ in (7), since we reduced modulo $Z$. Using (7), it 
would be enough to prove that 
$$ 
\deg \l((U\cap I^nV)/(U\cap I^nW))< \deg \l((U+I^nW)/I^nW).
$$
By the proof of Lemma 2(c), we know that 
$\deg \l((U\cap I^nV)/(U\cap I^nW))$
occurring in (7) and (8) can not exceed $\ell_{R/J_1}(I)-1$. Here we set 
$J_1:=ann_{R}(U\cap I^kV)$, so $J_1\supseteqq J$. Then
$$
\deg \l((U\cap I^nV)/(U\cap I^nW))\leq \ell_{R/J_1}(I)-1\leq 
\ell_{R/J}(I)-1=\ell_{N}(I)-1<\dim (H_1(\Cal{C}\otimes N)).
$$
But $\dim (H_1(\Cal{C}\otimes N))=\dim U=\deg \l((U+I^nW)/I^nW)$,
according to the proof of 2(c), and the claim is proven.
Similarly, $\deg \l(I^nV/I^nW)< \deg \l((U+I^nW)/I^nW)$. 
Now (7) shows that
$\deg \l(H_1(\Cal{C}\otimes \frac{N}{I^nN})) = \dim (H_1(\Cal{C}\otimes N))$. 
\enddemo

\subhead{Remark}
\endsubhead
The following Corollary generalizes Kodiyalam's results for $\t$ and $\e$.

\proclaim{Corollary 4} Let $(R,m)$ be Noetherian local, $I$ an ideal and $N$, 
$M$ finitely generated $R$-modules.
\roster
\item"(a)" If $\t_i (N/I^nN,M)$ has finite length for some $i$ 
and $n\gg 0$, then $\l(\t_i (N/I^nN,M))$ has polynomial growth for 
all large $n$. Moreover, 
$$
\deg \l(\t_i (N/I^nN,M))\leq \max\{\dim \t_i (N,M), \ell_{N}(I)-1\}.
$$ 
Equality holds, if $\dim \t_i (N,M)\geq \ell_{N}(I)$.

\item"(b)" If $\e^i(M,N/I^nN)$ has finite length for some $i$ and 
$n\gg 0$, then  $\l(\e^i(M,N/I^nN))$ has polynomial growth for all 
large $n$. Moreover, 
$$
\deg \l(\e^i(M,N/I^nN))\leq \max\{\dim \e^i (M,N), \ell_{N}(I)-1\}.
$$ 
Equality holds, if $\dim \e^i (M,N)\geq \ell_{N}(I)$.
\endroster
\noindent 
In particular, if $Supp(N)\cap Supp(M)\cap V(I)
=\{m\}$, then both $\l(\t_i (N/I^nN,M))$ and
$\l(\e^i(M,N/I^nN))$ have polynomial growth.
\endproclaim

\demo{Proof}\newline
(a) Apply 3(a) and (c) to $\Cal C$, a free resolution of $M$, consisting of 
finitely generated modules.\newline
(b) If $\tilde{\Cal C}$ is a free resolution of $M$, 
$\e_R^i (M,N/I^nN)$ is the $i^{-th}$ cohomology module of the complex
$$
\h (\tilde{\Cal C},\frac{N}{I^nN}) \cong  \h (\tilde{\Cal C},R)\otimes 
\frac{N}{I^nN} \cong (\h (\tilde{\Cal C},R) \otimes N) \otimes \frac{R}{I^n}
\cong (\h (\tilde{\Cal C},N) \otimes \frac{R}{I^n}.
$$
\noindent
Since $\Cal C:= \h (\tilde{\Cal C},N)$ is a complex of finitely generated 
modules, Proposition 3 (a) and (c) applies. 
\enddemo
\subhead{Example}
\endsubhead
(a) In Corollary 4(a) above, take $i=1$, $I=m$, $M=R/m$ and $N=R$.
Then $\t_1(R/m,R/m^n)=m^n/m^{n+1}$, which shows that 
$\deg (\l (\t_1(R/m,R/m^n))) = \deg (\l (m^n/m^{n+1})) = \ell (I)-1$.\newline
\noindent
(b) Let $(A,\tilde{m})$ be local, and $J_1$, $J_2$ be two 
$\tilde{m}$-primary ideals of $A$, such that $J_1J_2\subsetneq J_1\cap J_2$.
Let $X$ be an indeterminate over $A$. Define $R:= A[X]_{(\tilde{m},X)}$, 
and take in Corollary 4(a) $i=1$, $M=R/J_1R$, $N=R/J_2R$, $I=XR$. Then
$$
\align \t_1(N/I^nN,M)&=\t_1((R/J_2R)/X^n(R/J_2R),R/J_1R)\\
=\t_1(R/(X^nR+J_2R),R/J_1R)&=J_1R\cap (X^nR+J_2R)/J_1(X^nR+J_2R).
\endalign
$$
         
The numerator of the latter expression can be identified with all polynomials 
over $A$ having coefficients in $J_1\cap J_2$ in degree less than $n$, and
for which the coefficients in degree $n$ and higher are in $J_1$.
Similarly, the denominator can be viewed as polynomials over $A$, having 
coefficients in $J_1J_2$ in degree less than $n$, and all the other 
coefficients in $J_1$. Then the latter quotient above can be thought of as
polynomials of degree less than $n$, with coefficients in 
$(J_1\cap J_2)/J_1J_2$. In other words,
$\t_1(N/I^nN,M)$ is isomorphic to $n$ copies of $\t_1^A(A/J_2,A/J_1)$.
Therefore $\deg \l(\t_1(N/I^nN,M)=1$.
On the other hand, it is easy to see that $\dim \t_1(R/J_2R,R/J_1R)=1
>\ell_{R/J_2R}(XR)-1=0$. By Proposition 3(c) then, $\deg \l(\t_1(N/I^nN,M))=1$.
Hence, (a) and (b) show that both possible upper bounds given in 
4(a) are actually attained.

\head{The main result}
\endhead
In this section we show that, for finite $R$-modules $M$ and $N$,
$\l(\e^i(N/I^nN,M))$ is given by a polynomial for large $n$, 
provided
$$
  V(I) \cap Supp(M) \cap Supp(N)= \Cal{S},
$$
where $\Cal{S}$ is a finite subset of the set of maximal ideals of $R$.
This result is given in Corollary 6, while Theorem 5 gives the result 
in its maximum generality.\newline

 Recall that for an $R$-module $M$, $E_R(M)$ denotes the injective hull of 
$M$ over $R$. Also, $\mu^j(P,M)$ denotes the j-th Bass number of $M_P$,
that is, the number of copies of $E_R(R/P)$ occuring at the j-th place
in the minimal injective resolution of $M$. 
 The following Theorem is the most general form of our result:

\proclaim{Theorem 5} Let $R$ be a Noetherian ring, $I$ an ideal, $N$ a 
finitely generated $R$-module and $M$ an $R$-module such that 
$V(I) \cap Supp(M) \cap Supp(N)= \Cal S$, a finite subset of maximal 
ideals of $R$.
Assume that for all $ m \in \Cal S$, we have $\mu^j (m,M)<\infty$, 
where $j \in \{ i-1,i\}$ for some $i\geq 0$.
If $\e_R^i(N/I^nN,M)$ has finite length for large $n$, then 
$\l (\e_R^i (N/I^nN,M))$ coincides eventually with a polynomial.
\endproclaim
\demo{Proof}
Since $N$ is finitely generated, we can reduce the problem to the case 
$(R,m)$ is Noetherian and local, and $\Cal S=\{m\}$, by localizing at 
any prime $P$ in $\Cal S$ and using the formula

$$
\l (\e^i (\frac{N}{I^nN},M))= \sum_P \l (\e^i (\frac{N}{I^nN},M))_P=
\sum_P \l (\e_{R_P}^i (\frac{N_P}{I_P^nN_P},M_P).
$$
\noindent
Here the sum ranges over $V(I) \cap Supp(M) \cap Supp(N)$.

Let now $\tilde {\Cal C}$ be a minimal injective resolution of $M$. Then 
$Supp (\tilde {\Cal C}) \subset Supp(M)$, where \newline
$Supp (\tilde {\Cal C}) = \{ P \in Spec (R) : \tilde { \Cal C_P}$ is not the 
zero complex $\}.$ In other words, if $M_P=0$ for some $P \in Spec (R)$, then 
$\tilde {\Cal C_P}$ vanishes, too, by minimality of the injective resolution.
Now recall that every injective $R$-module is a direct sum of indecomposable 
injectives of the form $E_R(R/P)$, where $P\in Spec(R)$. It follows that
all prime ideals of $R$ associated to some module in the minimal injective 
resolution of $M$ are among those in $Supp(M)$. Since $N/I^nN$ is finitely 
generated,

$$
Supp (\h (\frac{N}{I^nN}, \tilde {\Cal C})) \subset Supp(N) \cap V(I) \cap 
Supp (\tilde{\Cal C}) \subset Supp(N) \cap V(I) \cap Supp(M) = \{ m \}.
$$
\noindent

A typical injective in $\tilde {\Cal C}$ has the form 
$E_l=\bigoplus_{P\in supp(M)}E_R(R/P)^{\mu^l(P,M)}$, $l\geq 0$, and since 
$\h(N/I^nN,-)$ distributes over direct sums, we get that 
the corresponding module in $\h(N/I^nN,\tilde {\Cal C} )$ is isomorphic to
$\bigoplus_{P\in Supp(M)}\h (N/I^nN,E_R(R/P))^{\mu^l(P,M)}$. But our support 
hypothesis says that every such module $\h (N/I^nN,E_R(R/P))$
vanishes, unless $P=m$.
Hence, by the above remarks on the supports, we see that

$$
\h (\frac{N}{I^nN},\tilde{\Cal C})=\h (\frac{N}{I^nN},\Cal C),
$$

\noindent
as complexes, where $\Cal C :=\Gamma_m ( \tilde{\Cal C}) $,    
with $\Gamma_m$ the $0^{-th}$ local cohomology functor. Thus, $\Cal C$
is the subcomplex of $\tilde {\Cal C}$ whose modules consist of
direct sums of $E_R(R/m)$. The finiteness assumption on the Bass numbers 
$\mu^i(m,M)$ and $\mu^{i-1}(m,M)$ simply says that there are only finitely 
many terms in $\Cal C$, at the $i^{-th}$ and $(i-1)^{-st}$ locations. 
Hence the two modules at the $i^{-th}$ and $(i-1)^{-st}$ locations in 
$\Cal C$ are Artinian modules.  
Finally, note that

$$
\h (\frac{N}{I^nN},\Cal C)= \h (\frac{R}{I^n}, \h (N, \Cal C) )
$$

\noindent
and $\h (N,\Cal C)$ is a complex with two Artinian modules at the 
$i^{-th}$ and $(i-1)^{-st}$ locations. 

Passing to the Matlis dual preserves the length of the (co)homology modules.
Thus,

$$
\l(\e^i(\frac{N}{I^nN},M)):=\l(H^i(\h(\frac{N}{I^nN},
\tilde{\Cal C})))=\l(H^i(\h(\frac{N}{I^nN},\Cal C)))=
\tag 10
$$
$$
\l(H^i(\h(\frac{R}{I^n}, \h (N,\Cal C))))
=\l(H_{i}(\h(\frac{R}{I^n},\h(N,\Cal C)))^{\c})=
\l(H_{i}(\h(N,\Cal C)^{\c}\otimes \frac{R}{I^n})),
$$ 

\noindent
the last equality by Lemma 1. 

The hypotheses in Proposition 3 are now met, 
for the two modules in $\h (N, \Cal C)^{\c}$ at the $i^{-th}$ and 
$(i-1)^{-st}$ locations are finitely generated over $\hat R$, the completion 
of $R$. Indeed, the Matlis dual of an Artinian $R$-module is finitely 
generated over $\hat R$. This concludes the proof of Theorem 5.
\enddemo

\proclaim{Corollary 6}Let $R$ be Noetherian, $I$ an ideal, $N$ a finitely 
generated $R$-module.
\roster
\item"(a)"
If $M$ is a finitely generated $R$-module such that $R/(I+annM+annN)$ is 
Artinian, then $\e^i (N/I^nN,M)$ has finite length for all\quad $i$ 
and  $n \geq 1$, and \newline 
$\l (\e^i (N/I^nN,M))$ has polynomial growth for all large $n$.
\item"(b)"
If $A$ is Artinian and $\e^i (N/I^nN,A)$ has finite length for some 
$i$ and $n \gg 0$, then  $\l (\e^i (N/I^nN,A))$ has
polynomial growth for all large $n$.
\item"(c)"
Assume $(R,m)$ is complete, and $M$ is Matlis reflexive (that is, 
$M \hookrightarrow ({M}^{\c})^{\c}$ is an isomorphism). 
Assume, furthermore, that $Supp(M) \cap Supp(N) \cap V(I) = \{m\}$.
If $\l (\e_R^i (N/I^nN,M))$ is finite for some $i$ and $n\gg 0$, then 
$\l (\e_R^i (N/I^nN,M))$ has polynomial growth for all large $n$.
\endroster    
\endproclaim

\demo{Proof}\newline
(a) Clearly, the condition $R/(I+ann M + ann N)$ Artinian is equivalent 
to the support condition given in Theorem 5, in the case where both $M$, $N$  
are finitely generated. Moreover, all the Bass numbers of $M$ are finite,
so Theorem 5 applies. \newline
(b) All the Bass numbers of an Artinian module $A$ are finite. 
To see this, recall that a module is Artinian if and only if it has finite 
socle and it is an essential extension of its socle. So, if $E$ is an 
injective hull of $A$, then it has the same socle as $A$ and, moreover, 
$E/A$ is still Artinian, since $E$ is so. It suffices to note that 
$\mu^0(m,A)=\mu^0(m,E)<\infty$, and $\mu^i(m,A)=\mu^{i-1}(m,E/A)$, 
for all $i\geq 1$.
Thus, the finiteness of Bass numbers for $A$ follows by induction, and 
Theorem 5 applies again.
Also note that (b) follows from 4(a) by the duality formula connecting
$\t $ and $\e $. \newline
(c) We may assume, without loss of generality, that $R$ is complete. 
It is known (see \cite {E}) that a module $M$ is Matlis reflexive if 
and only if there is an exact sequence

$$
   0@>>>K@>>>M@>>>A@>>>0,
$$
\noindent
where $A$ is an Artinian $R$-module and $K$ is a Noetherian $R$-module.
From the long exact sequence for $\e $ we see that all the Bass numbers of 
$M$ are finite, and Theorem 5 applies.
\enddemo

The following Corollary gives a degree estimate for the 
``harder'' $\e$.
\proclaim{Corollary 7}
Let $(R,m)$ be a local Noetherian ring of dimension $d$, $I$ an ideal, 
$M$, $N$  finite $R$-modules, such that $I+ann M+ ann N$ is $m$-primary. 
Let $\Cal C$ be a minimal injective resolution of $M$. Then, for each 
$i\geq 0$, we have 
$$
\deg\l(\e^i(N/I^nN,M))\leq\max\{\dim(H_i(\Gamma_m(\Cal{C})^{\c}\otimes N)), 
 \ell_{N}(I)-1\}.
$$
In particular, if $N=R$ and $0\leq i \leq \dim M$, then
$$
\deg \l(\e^i(R/I^n,M))\leq \max \{\dim R/annH_m^i(M), \ell(I)-1\}.
$$
\noindent
Moreover, if $R$ is Cohen-Macaulay, $i=d$ and $M=R$, then
$\deg \l(\e^d(R/I^n,R))=d$.
\endproclaim
\demo{Proof} The first statement follows from (10), using 
Proposition 3(c). The second statement follows then, since Matlis duality
preserves annihilators. For the third statement, note that we may assume that
$R$ is complete, and then $ann H_m^d(R)=ann (H_m^d(R))^{\c}=ann\, \omega_R$,
where $\omega_R$ denotes the canonical module of $R$. This has dimension 
$d$. On the other hand, $\ell (I)\leq d$, and Proposition 3(c) concludes
the proof.
\enddemo
\subhead{Remark}
\endsubhead
Kirby \cite {Kir} also showed in his paper that $\deg \l(\e^d(R/I^n,R))=d$, 
in the case where $R$ is Cohen-Macaulay and $I$ is $m$-primary.

Recall that, for an $R$-module $M$, $Soc(M):=(0:_Mm)$ and $Top(M):=M/mM$, 
where $m$ is the maximal ideal of a local ring $R$. In his paper, 
Kodiyalam also shows that, assuming $\l(M\otimes N)<\infty$, the top and 
the socle of $\t^i(N/I^nN,M)$ and $\e^i(M,N/I^nN)$ have lengths eventually 
given by polynomials. 
More can be shown in this direction.
With the technique already introduced, if $L$, $M$, $N$ are finite 
$R$-modules and $i$, $j$ are fixed nonegative numbers, the usual 
conclusion holds for each of the following: $\t^j(L,\t^i(N/I^nN,M))$ and 
$\e^j(L,\t^i(N/I^nN,M))$, $\t^j(L,\e^i(M,N/I^nN))$ and 
$\e^j(L,\e^i(M,N/I^nN))$. The only hypothesis needed is that they have 
finite length for all $n \gg 0$. This recovers Kodiyalam's top and socle 
results by taking $i=j=0$ and $L=k$, the residue field of $R$. 
In this paper we will only give the corresponding versions for 
$\e^i(N/I^nN,M)$, the module we are mainly interested in.
\proclaim{Corollary 8}Let $(R,m)$ be a local Noetherian ring, 
$L$, $M$, $N$ be finite $R$-modules, and $i$, $j$ fixed nonnegative 
integers. Assume $ann_RM+ann_RN+I$ is an $m$-primary ideal.
Then both $\t^j(L,\e^i(N/I^nN,M))$ and $\e^j(L,\e^i(N/I^nN,M))$ 
have finite lengths for all $n \geq 1$, given by polynomials for all 
large $n$.
\endproclaim

\demo{Proof}
We may again assume that $R$ is complete. Clearly, all the $\t$ and $\e$ 
modules in this Corollary have finite length for all $n$. By the proof 
of Theorem 5, we know that the Matlis dual of $\e^i(N/I^nN,M)$ equals 
$(U+I^nV)/I^nW$ for some finite $U$, $V$, $W$.
By applying the Matlis dual to both $\t^j(L,\e^i(N/I^nN,M))$ and 
$\e^j(L,\e^i(N/I^nN,M))$, we change the problem 
to showing that the conclusion of the Corollary holds for \newline 
$\e^j(L,(U+I^nV)/I^nW)$ and $\t^j(L,(U+I^nV)/I^nW)$, respectively. 
Let $\Cal C$ be a (minimal) free resolution of $L$. Then the $\e$ and $\t$ 
above are the $j^{-th}$ (co)homology modules of 
$\Cal C^T \otimes (U+I^nV)/I^nW$ and $\Cal C\otimes (U+I^nV)/I^nW$, 
respectively. Here $\Cal C^T$ stands for $\h(\Cal C,R)$. 
Proposition 3(b) now finishes the proof.
\enddemo

\head{An Example}
\endhead
The following example shows that $\e^i(N/I^nN,M)$, unlike $\e^i(M,N/I^nN)$,
may fail to have polynomial growth when we simply require that 
$\e^i(N/I^nN,M)$ have finite length for all large $n$.
Thus, in order to conclude that $\e^i(N/I^nN,M)$ has polynomial growth with 
respect to $n$, a more restrictive condition (such as the support condition 
in Theorem 5), beyond the obvious ``finite length'' condition, is necessary.
 
Consider $R=k[[X^2, XY, Y^2]]$, a subring of the power series
ring $k[[X,Y]]$, where $k$ is a field. It is easy to see that 
$R \cong k[[U, V, W]]/(V^2-UW)$, so $R$ is a two dimensional, Gorenstein, 
local ring. 
Take  $I=(X^2, XY)R$, which clearly is a height one prime of $R$, 
and note that the maximal ideal $m = (X^2, XY, Y^2)$ is associated 
to $I^n$ for all $n>1$.
Indeed, 
$$
I^n=(X^{2n},\quad X^{2n-1}Y, \quad X^{2n-2}Y^2, \quad X^{2n-3}Y^3, \hdots, 
X^nY^n)R
$$
\noindent
and for $a=X^{2n-2}$, $a \in R \smallsetminus I^n$, we have $am\subseteq I^n$. 
Finally, take $M=N=R$. By local duality, we have 
$$ 
\e^2(R/I^n,R)^{\c}= H_m^0(R/I^n)
$$
\noindent
and the latter has finite length for all $n \geq 2$.
Thus, $\e^2(R/I^n,R)$ has finite length for all $n \geq 2$, even though
$Supp (M)\cap Supp (N)\cap V(I)=\{I,m\}$ in this case.
On the other hand, 

$$
H_m^0(R/I^n)=H_{(Y^2)R}^0(R/I^n) \cong \frac{\cup_{k=1}^{\infty} 
(I^n \:_R \quad Y^{2k})} {I^n}
$$
\noindent
is generated by (classes of) monomials.
Note also that $H_m^0(R/I^n)$ is equal to $I^{(n)}/I^n$.
 
 Assume first that $n$ is odd.
Clearly, by dividing all the monomials generating $I^n$  by $Y^2$, we get all 
the monomials in $I^{(n)}\smallsetminus I^n$ that are multiplied into $I^n$ by 
$Y^2$, namely

$$
  X^{2n-2},\quad  X^{2n-3}Y, \hdots , X^{n+1}Y^{n-3},\quad  X^nY^{n-2}.  
$$ 
\noindent
By repeating the procedure on these monomials, we get those monomials 
in $I^{(n)}\smallsetminus I^n$ that are multiplied into $I^n$ by $Y^4$ :

$$
  X^{2n-4},\quad  X^{2n-5}Y, \hdots , X^{n+1}Y^{n-5},\quad  X^nY^{n-4}.
$$
\noindent
Continuing, we can get all the monomials generating $I^{(n)}/I^n$,          
the last ones being

$$X^{n+1},X^nY.$$ 
\noindent
If $n$ is even, we similarly get that $I^{(n)}/I^n$ is generated as a vector
space by
$$
  X^{2n-2},\quad  X^{2n-3}Y, \hdots , X^{n+1}Y^{n-3},\quad  X^nY^{n-2},  
$$
\noindent
$$
  X^{2n-4},\quad  X^{2n-5}Y, \hdots , X^{n+1}Y^{n-5},\quad  X^nY^{n-4},
$$
\noindent 
and so on, the last generator being now
$$X^n.$$  
\noindent
Note that the monomials listed in each case actually form a basis of 
$I^{(n)}/I^n$ as a k-vector space, so that the desired length is simply the 
number of these monomials.
Thus,

$$
\l(Ext_R^2(R/I^n,R))=2+4+6+ \hdots +(n-3)+(n-1)=\frac{n^2-1}{4}
$$
\noindent
for odd $n$, while
$$
\l(Ext_R^2(R/I^n,R))=1+3+5+ \hdots +(n-3)+(n-1)=\frac{n^2}{4}
$$
\noindent
for even $n$.
This not only says that $\l(Ext_R^2(R/I^n,R))$ is no longer a polynomial,
but also that the leading ``normalized'' coefficient (that gives multiplicity
in the classic theory) is no longer an integer. Here the normalized leading 
coefficient equals $\frac{1}
{2}$. It would be interesting to know if in situations such as the one in the
$Example$, the length formula is at least a ``periodic polynomial'', that is,
a polynomial with periodic coefficients.

\enddocument